# Some integrals and series involving the Stieltjes constants

Donal F. Connon

dconnon@btopenworld.com

16 January 2018

**Abstract**

Inter alia, we present a Fourier series involving the generalised Stieltjes constants.

## 1. Introduction

The Hurwitz zeta function $\varsigma(s,x)$ is initially defined for $\operatorname{Re} s > 1$ and $x > 0$ by

$$(1.1) \qquad \varsigma(s,x) = \sum_{n=0}^{\infty} \frac{1}{(n+x)^s}$$

and $\varsigma'(s,x) := \frac{\partial}{\partial s} \varsigma(s,x)$. Note that $\varsigma(s,x)$ may be analytically continued to the whole $s$ plane except for a simple pole at $s = 1$. For example, Hasse (1898-1979) showed that [19]

$$(1.2) \qquad \varsigma(s,x) = \frac{1}{s-1} \sum_{n=0}^{\infty} \frac{1}{n+1} \sum_{k=0}^{n} \binom{n}{k} \frac{(-1)^k}{(k+x)^{s-1}}$$

is a globally convergent series for $\varsigma(s,x)$ and, except for $s = 1$, provides an analytic continuation of $\varsigma(s,x)$ to the entire complex plane.

It may be noted from (1.1) that $\varsigma(s,1) = \varsigma(s)$.

We easily see from (1.2) that

$$\lim_{s \to 1}[(s-1)\varsigma(s,x)] = \sum_{n=0}^{\infty} \frac{1}{n+1} \sum_{k=0}^{n} \binom{n}{k} (-1)^k$$

and, since $(1-1)^n = \sum_{k=0}^{n} \binom{n}{k}(-1)^k = \delta_{n,0}$, we have $\lim_{s \to 1}[(s-1)\varsigma(s,x)] = 1$ (which shows that $\varsigma(s,x)$ has a simple pole at $s = 1$).

We may therefore form the Laurent expansion of the Hurwitz zeta function $\varsigma(s,x)$ about $s = 1$

$$(1.3) \qquad \varsigma(s,x) = \frac{1}{s-1} + \sum_{n=0}^{\infty} \frac{(-1)^n}{n!} \gamma_n(x)(s-1)^n$$

where $\gamma_n(x)$ are the known as the Stieltjes constants and $\gamma_n(1) \equiv \gamma_n$.

## 2. The generalised Stieltjes constants $\gamma_n(x)$

Using (1.1) we see that

(2.1) $$\frac{\partial}{\partial x}\varsigma(s,x) = -s\varsigma(s+1,x)$$

we see from (1.3) that

$$\frac{\partial}{\partial x}\varsigma(s,x) = -1 - \sum_{n=0}^{\infty}\frac{(-1)^n}{n!}\gamma_n(x)s^{n+1}$$

and hence we have a fundamental formula for $m \geq 0$

(2.2) $$\left.\frac{\partial^{m+1}}{\partial s^{m+1}}\frac{\partial}{\partial x}\varsigma(s,x)\right|_{s=0} = (-1)^{m+1}(m+1)\gamma_m(x)$$

Differentiating (1.2) gives us

$$\frac{\partial}{\partial x}\varsigma(s,x) = -\sum_{n=0}^{\infty}\frac{1}{n+1}\sum_{k=0}^{n}\binom{n}{k}\frac{(-1)^k}{(k+x)^s}$$

and

$$\frac{\partial^{m+1}}{\partial s^{m+1}}\frac{\partial}{\partial x}\varsigma(s,x) = (-1)^m\sum_{n=0}^{\infty}\frac{1}{n+1}\sum_{k=0}^{n}\binom{n}{k}(-1)^k\frac{\log^{m+1}(k+x)}{(k+x)^s}$$

We therefore obtain

(2.3) $$\gamma_m(x) = -\frac{1}{m+1}\sum_{n=0}^{\infty}\frac{1}{n+1}\sum_{k=0}^{n}\binom{n}{k}(-1)^k\log^{m+1}(k+x)$$

This is a much more succinct derivation of the formula originally reported in [14] in 2007.

□

Using (2.1) with (1.3) we obtain

$$\frac{\partial}{\partial x}\varsigma(s,x) = \sum_{n=0}^{\infty}\frac{(-1)^n}{n!}\gamma_n'(x)(s-1)^n$$

and using (2.1) we have

(2.4) $$s\varsigma(s+1,x) = -\sum_{n=0}^{\infty}\frac{(-1)^n}{n!}\gamma_n'(x)(s-1)^n$$

Evaluation at $s=1$ gives us

(2.5)
$$\varsigma(2,x) = -\gamma_0'(x)$$

Using the Weierstrass canonical form of the gamma function [27, p.1] we have

(2.6)
$$\log \Gamma(1+x) = -\gamma x - \sum_{n=1}^{\infty}\left[\log\left(1+\frac{x}{n}\right)-\frac{x}{n}\right]$$

and differentiation gives us

(2.7)
$$\psi(1+x) = -\gamma + \sum_{n=1}^{\infty}\left(\frac{1}{n}-\frac{1}{n+x}\right)$$

where $\psi(u)$ is the digamma function which is the logarithmic derivative of the gamma function $\psi(u) := \frac{d}{du}\log\Gamma(u)$.

From (2.7) we may easily deduce the functional equation

(2.8)
$$\psi(1+x) = \psi(x) + \frac{1}{x}$$

and differentiation of (2.7) gives us

$$\psi^{(k)}(x) = (-1)^{k+1} k! \varsigma(k+1,x)$$

In particular we have

$$\varsigma(2,x) = \psi'(x)$$

and (2.5) gives us

$$\psi'(x) = -\gamma_0'(x)$$

Integration results in

$$\gamma_0(x) = -\psi(x) + c$$

To determine the constant of integration we use the elementary identity (which may be easily deduced from the definition in (1.1))

$$\varsigma(s,x) = \varsigma(s,1+x) + \frac{1}{x^s}$$

from which we deduce that

$$\lim_{s \to 1}\left[\varsigma(s,x) - \varsigma(s,1+x)\right] = \frac{1}{x}$$

and, using (1.3), we see that

$$\gamma_0(x) - \gamma_0(1+x) = \lim_{s \to 1}[\varsigma(s,x) - \varsigma(s,1+x)]$$

Hence, we have

(2.9) $$\gamma_0(x) - \gamma_0(1+x) = \frac{1}{x}$$

This then shows us that $c = 0$ and we therefore obtain the well-known formula [30]

(2.10) $$\gamma_0(x) = -\psi(x)$$

Using (2.3) we see that

(2.11) $$\psi(x) = \sum_{n=0}^{\infty} \frac{1}{n+1} \sum_{k=0}^{n} \binom{n}{k} (-1)^k \log(k+x)$$

as first noted by Sondow and Hadjicostas [26].

$\square$

Coffey [11] has reported that for $n \geq 1$

$$\int_0^1 \frac{u^{x-1}}{\log u}(1-u)^n \, du = \sum_{k=0}^{n} \binom{n}{k} (-1)^k \log(k+x)$$

and we note that the integrand is negative in the interval $(0,1]$.

We have the formal summation

$$\sum_{n=1}^{\infty} \frac{1}{n+1} \int_0^1 \frac{u^{x-1}}{\log u}(1-u)^n \, du = \sum_{m=2}^{\infty} \frac{1}{m} \int_0^1 \frac{u^{x-1}}{(1-u)\log u}(1-u)^m \, du$$

$$= -\int_0^1 \frac{u^{x-1}}{(1-u)\log u}[\log u + (1-u)] \, du$$

$$= \int_0^1 u^{x-1} \left[\frac{1}{\log u} - \frac{1}{1-u}\right] du$$

We see from (2.3) that

$$\psi(x) = \log x + \sum_{n=1}^{\infty} \frac{1}{n+1} \sum_{k=0}^{n} \binom{n}{k} (-1)^k \log(k+x)$$

and hence we obtain the known integral

$$\psi(x) - \log x = \int_0^1 u^{x-1} \left[ \frac{1}{\log u} - \frac{1}{1-u} \right] du$$

In addition, we also see that $0 > \psi(x) - \log x$ for all $x > 0$.

□

Differentiating (2.4) with respect to $s$ gives us

$$s\varsigma'(s+1, x) + \varsigma(s+1, x) = -\sum_{n=0}^{\infty} \frac{(-1)^n}{n!} \gamma_n'(x) n(s-1)^{n-1}$$

and we obtain with $s = 1$

$$\gamma_1'(x) = \varsigma'(2, x) + \varsigma(2, x)$$

as was previously derived in a different manner in [15].

Therefore, we have

$$\gamma_1'(x) = \sum_{k=0}^{\infty} \frac{1 - \log(k+x)}{(k+x)^2}$$

and we may deduce that

$$\gamma_1'(x) < 0 \text{ for all } x \geq e$$

We therefore determine that $\gamma_1(x)$ is monotonic decreasing for (at least) all $x \geq e$. However, we note that $\gamma_1(x)$ is not monotonic decreasing throughout $(0,1]$ because, for example, we have [12]

$$\gamma_1\left(\frac{1}{2}\right) = \gamma_1 - \log^2 2 - 2\gamma \log 2$$

**3. A Fourier series involving the generalised Stieltjes constants**

Lerch [23] proved the following theorem in 1895:

If the series

(3.1) $$f(x) = \sum_{n=1}^{\infty} \frac{c_n}{n} \sin 2\pi nx$$

is convergent for $0 < x < 1$, then the derivative of $f(x)$ is given in this interval by

(3.2) $$f'(x) \cdot \frac{\sin \pi x}{\pi} = \sum_{n=0}^{\infty} (c_n - c_{n+1}) \sin(2n+1)\pi x$$

where $c_0 = 0$, provided that the last series converges uniformly for $\varepsilon \leq x \leq 1-\varepsilon$ for all $\varepsilon > 0$. Equation (3.2) may be written as

$$(3.3) \qquad f'(x) \cdot \frac{\sin \pi x}{\pi} = -c_1 \sin \pi x + \sum_{n=1}^{\infty} (c_n - c_{n+1}) \sin(2n+1)\pi x$$

Subject to the same conditions, if the series

$$(3.4) \qquad g(x) = \sum_{n=1}^{\infty} \frac{d_n}{n} \cos 2\pi n x$$

is convergent for $0 < x < 1$, then we have

$$(3.5) \qquad g'(x) \cdot \frac{\sin \pi x}{\pi} = \sum_{n=0}^{\infty} (d_n - d_{n+1}) \cos(2n+1)\pi x$$

where $d_0 = 0$. Or equivalently

$$(3.6) \qquad g'(x) \cdot \frac{\sin \pi x}{\pi} = -d_1 \cos \pi x + \sum_{n=1}^{\infty} (d_n - d_{n+1}) \cos(2n+1)\pi x$$

Suppose $f(x) + g(x)$ may be represented by

$$f(x) + g(x) = \sum_{n=1}^{\infty} \frac{c_n}{n} \sin 2\pi n x + \sum_{n=1}^{\infty} \frac{d_n}{n} \cos 2\pi n x$$

and assume that there exists a function $\phi(x)$ such that

$$f(x) + \phi(x) = \sum_{n=1}^{\infty} \frac{c_n}{n} \sin 2\pi n x$$

Then we have

$$g(x) - \phi(x) = \sum_{n=1}^{\infty} \frac{d_n}{n} \sin 2\pi n x$$

and hence we see that

$$(3.7) \qquad [f'(x) + g'(x)] \cdot \frac{\sin \pi x}{\pi} = -c_1 \sin \pi x - d_1 \cos \pi x$$

$$+ \sum_{n=1}^{\infty} (c_n - c_{n+1}) \sin(2n+1)\pi x + \sum_{n=1}^{\infty} (d_n - d_{n+1}) \cos(2n+1)\pi x$$

$\square$

For example, letting $c_n \equiv 1$ in (3.1) gives us

$$f(x) = \sum_{n=1}^{\infty} \frac{\sin 2\pi n x}{n}$$

and employing (3.3) we have

$$f'(x) \cdot \frac{\sin \pi x}{\pi} = -\sin \pi x$$

Integration gives us

$$f(x) = -\pi x + c$$

and since $f(1/2) = 0$ we deduce that

$$f(x) = \frac{1}{2}\pi(1-2x)$$

Hence, we have the well-known Fourier series valid for $0 < x < 1$ as shown in Carslaw's book [8, p.241]

(3.8) $$\pi\left(\frac{1}{2} - x\right) = \sum_{n=1}^{\infty} \frac{\sin 2n\pi x}{n}$$

□

Similarly, letting $d_n \equiv 1$ in (3.4) results in

$$g(x) = \sum_{n=1}^{\infty} \frac{\cos 2\pi n x}{n}$$

and employing (3.6) we have

$$g'(x) \cdot \frac{\sin \pi x}{\pi} = -\cos \pi x$$

Integration gives us

$$g(x) = -\log \sin \pi x + c$$

and since $g(1/2) = \sum_{n=1}^{\infty} \frac{(-1)^n}{n} = -\log 2$ we deduce that

$$g(x) = -\log(2 \sin \pi x)$$

Hence, we have [8, p.241] for $0 < x < 1$

$$(3.9) \qquad -\log(2\sin \pi x) = \sum_{n=1}^{\infty} \frac{\cos 2n\pi x}{n}$$

These well-known Fourier series were derived in this manner by Lerch [23].

□

We have the well-known Hurwitz's formula for the Fourier expansion of the Hurwitz zeta function $\varsigma(s,x)$ as reported in Titchmarsh's treatise [28, p.37]

$$(3.10) \qquad \varsigma(s,x) = 2\Gamma(1-s)\left[\sin\left(\frac{\pi s}{2}\right)\sum_{n=1}^{\infty}\frac{\cos 2n\pi x}{(2\pi n)^{1-s}} + \cos\left(\frac{\pi s}{2}\right)\sum_{n=1}^{\infty}\frac{\sin 2n\pi x}{(2\pi n)^{1-s}}\right]$$

where $\operatorname{Re}(s) < 0$ and $0 < x \leq 1$. In 2000, Boudjelkha [6] showed that this formula also applies in the region $\operatorname{Re}(s) < 1$. It may be noted that when $x = 1$ this reduces to Riemann's functional equation for $\varsigma(s)$. Letting $s \to 1-s$ we may write (3.10) as

$$(3.11) \qquad \varsigma(1-s,x) = 2\Gamma(s)\left[\cos\left(\frac{\pi s}{2}\right)\sum_{n=1}^{\infty}\frac{\cos 2n\pi x}{(2\pi n)^{s}} + \sin\left(\frac{\pi s}{2}\right)\sum_{n=1}^{\infty}\frac{\sin 2n\pi x}{(2\pi n)^{s}}\right]$$

I once thought of using (3.10) to develop a Fourier series for the generalised Stieltjes constants but, unfortunately, (2.2) cannot be applied indiscriminately to Hurwitz's formula because this results in a divergent series. This explains Ramanujan's erroneous thinking in this area as exemplified in Berndt's book, *Ramanujan's Notebooks*, Part I [3, p.200].

Differentiation of (3.10) gives us

$$(3.12) \quad \varsigma'(s,x) = 2\Gamma(1-s)\left[\sin\left(\frac{\pi s}{2}\right)\sum_{n=1}^{\infty}\frac{\log(2\pi n)\cos 2n\pi x}{(2\pi n)^{1-s}} + \cos\left(\frac{\pi s}{2}\right)\sum_{n=1}^{\infty}\frac{\log(2\pi n)\sin 2n\pi x}{(2\pi n)^{1-s}}\right]$$

$$+\pi\Gamma(1-s)\left[\cos\left(\frac{\pi s}{2}\right)\sum_{n=1}^{\infty}\frac{\cos 2n\pi x}{(2\pi n)^{1-s}} - \sin\left(\frac{\pi s}{2}\right)\sum_{n=1}^{\infty}\frac{\sin 2n\pi x}{(2\pi n)^{1-s}}\right]$$

$$-2\Gamma'(1-s)\left[\sin\left(\frac{\pi s}{2}\right)\sum_{n=1}^{\infty}\frac{\cos 2n\pi x}{(2\pi n)^{1-s}} + \cos\left(\frac{\pi s}{2}\right)\sum_{n=1}^{\infty}\frac{\sin 2n\pi x}{(2\pi n)^{1-s}}\right]$$

Letting $s = 0$ gives us

$$(3.13) \qquad \varsigma'(0,x) = \frac{1}{\pi}\sum_{n=1}^{\infty}\frac{\log(2\pi n)+\gamma}{n}\sin 2n\pi x + \frac{1}{2}\sum_{n=1}^{\infty}\frac{\cos 2n\pi x}{n}$$

where we have used $\Gamma'(1) = \gamma$.

Using (3.8) and (3.9) and Lerch's identity [2]

(3.14) $$\log \Gamma(x) = \varsigma'(0,x) + \frac{1}{2}\log(2\pi)$$

we obtain Kummer's Fourier series [29] for the log gamma function

(3.15) $$\log \Gamma(x) = \frac{1}{2}\log\frac{\pi}{\sin \pi x} + [\gamma + \log(2\pi)]\left(\frac{1}{2} - x\right) + \frac{1}{\pi}\sum_{n=1}^{\infty}\frac{\log n}{n}\sin 2\pi n x$$

which, because we relied on (3.8) and (3.9), is only valid for $0 < x < 1$. As recently noted by Blagouchine [5], Kummer's Fourier series was actually first discovered by Malmstén [24] in 1846.

This may be written as

$$\log \Gamma(x) + \frac{1}{2}\log \sin \pi x - \frac{1}{2}\log \pi - [\gamma + \log(2\pi)]\left(\frac{1}{2} - x\right) = \frac{1}{\pi}\sum_{n=1}^{\infty}\frac{\log n}{n}\sin 2\pi n x$$

and applying Lerch's formula (3.3) gives us

(3.16) $$\sum_{n=1}^{\infty}\log\left(1+\frac{1}{n}\right)\sin(2n+1)\pi x = -\left[\psi(x)\sin \pi x + \frac{\pi}{2}\cos \pi x + (\gamma + \log 2\pi)\sin \pi x\right]$$

which is valid for $0 < x < 1$. See for example [18, p.105], [23] and [25, p.204].

Letting $x = \frac{1}{2}$ results in

(3.17) $$\log \frac{\pi}{2} = \sum_{n=1}^{\infty}(-1)^{n+1}\log\left(1+\frac{1}{n}\right)$$

being the logarithm of Wallis's product formula for $\pi/2$ which was published in Arithmetica Infinitorum in 1659

$$\frac{\pi}{2} = \prod_{n=1}^{\infty}\frac{(2n)^2}{(2n-1)(2n+1)}$$

□

The second derivative of (3.10) gives us

(3.18)
$$\varsigma''(s,x) = 2\Gamma(1-s)\left[\sin\left(\frac{\pi s}{2}\right)\sum_{n=1}^{\infty}\frac{\log^2(2\pi n)\cos 2n\pi x}{(2\pi n)^{1-s}} + \cos\left(\frac{\pi s}{2}\right)\sum_{n=1}^{\infty}\frac{\log^2(2\pi n)\sin 2n\pi x}{(2\pi n)^{1-s}}\right]$$

$$+\pi\Gamma(1-s)\left[\cos\left(\frac{\pi s}{2}\right)\sum_{n=1}^{\infty}\frac{\log(2\pi n)\cos 2n\pi x}{(2\pi n)^{1-s}} - \sin\left(\frac{\pi s}{2}\right)\sum_{n=1}^{\infty}\frac{\log(2\pi n)\sin 2n\pi x}{(2\pi n)^{1-s}}\right]$$

$$-2\Gamma'(1-s)\left[\sin\left(\frac{\pi s}{2}\right)\sum_{n=1}^{\infty}\frac{\log(2\pi n)\cos 2n\pi x}{(2\pi n)^{1-s}}+\cos\left(\frac{\pi s}{2}\right)\sum_{n=1}^{\infty}\frac{\log(2\pi n)\sin 2n\pi x}{(2\pi n)^{1-s}}\right]$$

$$+\pi\,\Gamma(1-s)\left[\cos\left(\frac{\pi s}{2}\right)\sum_{n=1}^{\infty}\frac{\log(2\pi n)\cos 2n\pi x}{(2\pi n)^{1-s}}-\sin\left(\frac{\pi s}{2}\right)\sum_{n=1}^{\infty}\frac{\log(2\pi n)\sin 2n\pi x}{(2\pi n)^{1-s}}\right]$$

$$-\frac{1}{2}\pi^2\,\Gamma(1-s)\left[\sin\left(\frac{\pi s}{2}\right)\sum_{n=1}^{\infty}\frac{\cos 2n\pi x}{(2\pi n)^{1-s}}+\cos\left(\frac{\pi s}{2}\right)\sum_{n=1}^{\infty}\frac{\sin 2n\pi x}{(2\pi n)^{1-s}}\right]$$

$$-\pi\,\Gamma'(1-s)\left[\cos\left(\frac{\pi s}{2}\right)\sum_{n=1}^{\infty}\frac{\cos 2n\pi x}{(2\pi n)^{1-s}}-\sin\left(\frac{\pi s}{2}\right)\sum_{n=1}^{\infty}\frac{\sin 2n\pi x}{(2\pi n)^{1-s}}\right]$$

$$-2\Gamma'(1-s)\left[\sin\left(\frac{\pi s}{2}\right)\sum_{n=1}^{\infty}\frac{\log(2\pi n)\cos 2n\pi x}{(2\pi n)^{1-s}}+\cos\left(\frac{\pi s}{2}\right)\sum_{n=1}^{\infty}\frac{\log(2\pi n)\sin 2n\pi x}{(2\pi n)^{1-s}}\right]$$

$$-\pi\,\Gamma'(1-s)\left[\cos\left(\frac{\pi s}{2}\right)\sum_{n=1}^{\infty}\frac{\cos 2n\pi x}{(2\pi n)^{1-s}}-\sin\left(\frac{\pi s}{2}\right)\sum_{n=1}^{\infty}\frac{\sin 2n\pi x}{(2\pi n)^{1-s}}\right]$$

$$+2\Gamma''(1-s)\left[\sin\left(\frac{\pi s}{2}\right)\sum_{n=1}^{\infty}\frac{\cos 2n\pi x}{(2\pi n)^{1-s}}+\cos\left(\frac{\pi s}{2}\right)\sum_{n=1}^{\infty}\frac{\sin 2n\pi x}{(2\pi n)^{1-s}}\right]$$

Letting $s=0$ gives us

$$(3.19)\ \varsigma''(0,x)=2\sum_{n=1}^{\infty}\frac{\log^2(2\pi n)\sin 2n\pi x}{2\pi n}+4\gamma\sum_{n=1}^{\infty}\frac{\log(2\pi n)\sin 2n\pi x}{2\pi n}+[2\gamma^2-\varsigma(2)]\sum_{n=1}^{\infty}\frac{\sin 2n\pi x}{2\pi n}$$

$$+2\pi\sum_{n=1}^{\infty}\frac{\log(2\pi n)\cos 2n\pi x}{2\pi n}+2\pi\,\gamma\sum_{n=1}^{\infty}\frac{\cos 2n\pi x}{2\pi n}$$

where we have used [27, p.265] $\Gamma''(1)=\gamma^2+\varsigma(2)$.

Letting $x\to 1-x$ in (3.19) results in

$$\varsigma''(0,1-x)=-2\sum_{n=1}^{\infty}\frac{\log^2(2\pi n)\sin 2n\pi x}{2\pi n}-4\gamma\sum_{n=1}^{\infty}\frac{\log(2\pi n)\sin 2n\pi x}{2\pi n}-[2\gamma^2-\varsigma(2)]\sum_{n=1}^{\infty}\frac{\sin 2n\pi x}{2\pi n}$$

$$+2\pi\sum_{n=1}^{\infty}\frac{\log(2\pi n)\cos 2n\pi x}{2\pi n}+2\pi\,\gamma\sum_{n=1}^{\infty}\frac{\cos 2n\pi x}{2\pi n}$$

In fact, without relying on (3.7), by using the representations

$$(3.20)\quad \varsigma(s,x)+\varsigma(s,1-x)=4\Gamma(1-s)\sin\left(\frac{\pi s}{2}\right)\sum_{n=1}^{\infty}\frac{\cos 2n\pi x}{(2\pi n)^{1-s}}$$

and

(3.21) $\zeta(s,x) - \zeta(s,1-x) = 4\Gamma(1-s)\cos\left(\frac{\pi s}{2}\right)\sum_{n=1}^{\infty}\frac{\sin 2n\pi x}{(2\pi n)^{1-s}}$

we can derive the equivalent forms

$$\zeta''(0,x) + \zeta''(0,1-x) = 2\sum_{n=1}^{\infty}\frac{\log(2\pi n) + \gamma}{n}\cos 2n\pi x$$

and

$$\zeta''(0,x) - \zeta''(0,1-x) = \sum_{n=1}^{\infty}\frac{2\log^2(2\pi n) + 4\gamma\log(2\pi n) + \left[2\gamma^2 - \zeta(2)\right]}{\pi n}\sin 2n\pi x$$

Deninger [17] showed that

(3.22) $\sum_{n=1}^{\infty}\frac{\log n}{n}\cos 2n\pi x = \frac{1}{2}\left[\zeta''(0,x) + \zeta''(0,1-x)\right] + \left[\gamma + \log(2\pi)\right]\log(2\sin\pi x)$

This is a companion to Kummer's Fourier series formula for $\log\Gamma(x)$ and using Euler's reflection formula $\Gamma(x)\Gamma(1-x) = \frac{\pi}{\sin\pi x}$ we may write this as

$$\sum_{n=1}^{\infty}\frac{\log n}{n}\sin 2n\pi x = \frac{\pi}{2}[\log\Gamma(x) - \log\Gamma(1-x)] - [\gamma + \log(2\pi)]\frac{\pi}{2}(1-2x)$$

With Lerch's identity (3.14) this becomes

$$\sum_{n=1}^{\infty}\frac{\log n}{n}\sin 2n\pi x = \frac{\pi}{2}[\zeta'(0,x) - \zeta'(0,1-x)] - [\gamma + \log(2\pi)]\frac{\pi}{2}(1-2x)$$

and, in this format, we can see the structural similarity with (5.7). Using (5.11) and (5.12) we may write

$$\sum_{n=1}^{\infty}\frac{\gamma + \log(2\pi n)}{n}\cos 2n\pi x = \frac{1}{2}\left[\zeta''(0,x) + \zeta''(0,1-x)\right]$$

$$\sum_{n=1}^{\infty}\frac{\gamma + \log(2\pi n)}{n}\sin 2n\pi x = \frac{\pi}{2}[\zeta'(0,x) - \zeta'(0,1-x)]$$

□

Applying Lerch's formula (3.7) to (3.19) we obtain

$$\frac{\sin\pi x}{\pi}\frac{d}{dx}\zeta''(0,x) = -c_1\sin\pi x + \sum_{n=1}^{\infty}(c_n - c_{n+1})\sin(2n+1)\pi x$$

$$-d_1 \cos \pi x + \sum_{n=1}^{\infty}(d_n - d_{n+1})\cos(2n+1)\pi x$$

where

$$c_n = \frac{1}{\pi}\log^2(2\pi n) + \frac{2\gamma}{\pi}\log(2\pi n) + \frac{1}{2\pi}[2\gamma^2 - \varsigma(2)]$$

$$d_n = \log(2\pi n) + \gamma$$

Using (2.2) this becomes

(3.23)
$$2\gamma_1(x)\frac{\sin \pi x}{\pi} = -c_1 \sin \pi x + \sum_{n=1}^{\infty}(c_n - c_{n+1})\sin(2n+1)\pi x$$

$$-d_1 \cos \pi x + \sum_{n=1}^{\infty}(d_n - d_{n+1})\cos(2n+1)\pi x$$

or equivalently

(3.24)
$$2\gamma_1(x)\frac{\sin \pi x}{\pi} = -c_1 \sin \pi x + \sum_{n=1}^{\infty}(c_n - c_{n+1})\sin(2n+1)\pi x$$

$$-[\log 2\pi + \gamma]\cos \pi x - \sum_{n=1}^{\infty}\log\left(1+\frac{1}{n}\right)\cos(2n+1)\pi x$$

For example, we have

$$2\gamma_1\left(\frac{1}{2}\right) = -\left(\log^2(2\pi) + 2\gamma\log(2\pi) + \frac{1}{2}[2\gamma^2 - \varsigma(2)]\right) + \pi\sum_{n=1}^{\infty}(-1)^n(c_n - c_{n+1})$$

and using Coffey's result [11]

$$\gamma_1\left(\frac{1}{2}\right) = \gamma_1 - \log^2 2 - 2\gamma\log 2$$

we may obtain a series representation involving $\gamma_1$.

It is known that ([5], [11], [15])

$$\gamma_1\left(\frac{p}{q}\right) = \gamma_1 - [\gamma + \log(2\pi)]\log(2\pi q) - \frac{1}{2}\log^2 q$$

$$+ \sum_{r=1}^{q-1}\varsigma''\left(0,\frac{r}{q}\right)\cos\left(\frac{2\pi rp}{q}\right) - 2[\gamma + \log(2\pi q)]\sum_{r=1}^{q-1}\log\Gamma\left(\frac{r}{q}\right)\cos\left(\frac{2\pi rp}{q}\right)$$

$$+\pi\sum_{r=1}^{q-1}\log\Gamma\left(\frac{r}{q}\right)\sin\left(\frac{2\pi rp}{q}\right)+\frac{\pi}{2}[\gamma+\log(2\pi q)]\cot\left(\frac{p\pi}{q}\right)$$

Some examples of the first Stieltjes constants are illustrated below:

$$\gamma_1\left(\frac{1}{4}\right)=\frac{1}{2}[2\gamma_1-7\log^2 2-6\gamma\log 2]-\frac{1}{2}\pi\left[\gamma+4\log 2+3\log\pi-4\log\Gamma\left(\frac{1}{4}\right)\right]$$

$$\gamma_1\left(\frac{1}{5}\right)=\frac{1}{4}\left[4\gamma_1-\frac{5}{2}\log^2 5-5\gamma\log 5\right]-\frac{1}{2}\pi[\log(10\pi)+\gamma]\cot\left(\frac{\pi}{5}\right)$$

$$+\frac{1}{4}\sqrt{5}\left[\varsigma''\left(0,\frac{1}{5}\right)-\varsigma''\left(0,\frac{2}{5}\right)-\varsigma''\left(0,\frac{3}{5}\right)+\varsigma''\left(0,\frac{4}{5}\right)-[\gamma+\log(2\pi)]\log\frac{1}{2}(3+\sqrt{5})\right]$$

Letting $x \to 1-x$ in (3.23) we see that

$$2\gamma_1(1-x)\frac{\sin\pi x}{\pi}=-c_1\sin\pi x+\sum_{n=1}^{\infty}(c_n-c_{n+1})\sin(2n+1)\pi x$$

$$+d_1\cos\pi x-\sum_{n=1}^{\infty}(d_n-d_{n+1})\cos(2n+1)\pi x$$

and thus

$$[\gamma_1(1-x)-\gamma_1(x)]\frac{\sin\pi x}{\pi}=d_1\cos\pi x-\sum_{n=1}^{\infty}(d_n-d_{n+1})\cos(2n+1)\pi x$$

Using $d_n=\log(2\pi n)+\gamma$ this may be written as

(3.25) $$\sum_{n=1}^{\infty}\log\left(1+\frac{1}{n}\right)\cos(2n+1)\pi x=[\gamma_1(1-x)-\gamma_1(x)]\frac{\sin\pi x}{\pi}-[\gamma+\log(2\pi)]\cos\pi x$$

Adamchik [1] has reported that

(3.36) $$\gamma_1\left(1-\frac{p}{q}\right)-\gamma_1\left(\frac{p}{q}\right)=\pi[\log(2\pi q)+\gamma]\cot\left(\frac{\pi p}{q}\right)-2\pi\sum_{j=1}^{q-1}\log\Gamma\left(\frac{j}{q}\right)\sin\left(\frac{2\pi jp}{q}\right)$$

where $p$ and $q$ are positive integers and $p<q$. Adamchik [1] notes that this formula was first proved by Almkvist and Meurman. Adamchik's derivation [1] is rather terse and an expanded exposition may be found in [15]. As recently noted by Blagouchine [5] this was in fact originally discovered by Carl Malmstén in 1846.

Hence, we obtain

(3.27)
$$\sum_{n=1}^{\infty} \log\left(1+\frac{1}{n}\right) \cos\left(\frac{(2n+1)\pi p}{q}\right) = \log q \cos\left(\frac{\pi p}{q}\right) - 2\sin\left(\frac{\pi p}{q}\right) \sum_{j=1}^{q-1} \log \Gamma\left(\frac{j}{q}\right) \sin\left(\frac{2\pi j p}{q}\right)$$

The above procedures may obviously be extended to $\gamma_m(x)$ but the algebra becomes increasingly more tedious.

To recap, we have

$$\sum_{n=1}^{\infty} \log\left(1+\frac{1}{n}\right) \cos(2n+1)\pi x = [\gamma_1(1-x) - \gamma_1(x)]\frac{\sin \pi x}{\pi} - [\gamma + \log(2\pi)]\cos \pi x$$

$$\sum_{n=1}^{\infty} \log\left(1+\frac{1}{n}\right) \sin(2n+1)\pi x = -\left[\psi(x)\sin \pi x + \frac{\pi}{2}\cos \pi x + [\gamma + \log 2\pi]\sin \pi x\right]$$

and solving a pair of simultaneous equations we obtain

(3.28) $$\sum_{n=1}^{\infty} \log\left(1+\frac{1}{n}\right) \cos 2n\pi x = [\gamma_1(1-x) - \gamma_1(x)]\frac{\sin \pi x \cos \pi x}{\pi} - [\gamma + \log(2\pi)]$$

$$-\left[\psi(x)\sin \pi x + \frac{\pi}{2}\cos \pi x\right]\sin \pi x$$

(3.29)
$$\sum_{n=1}^{\infty} \log\left(1+\frac{1}{n}\right) \sin 2n\pi x = -[\gamma_1(1-x) - \gamma_1(x)]\frac{\sin^2 \pi x}{\pi} - \left[\psi(x)\sin \pi x + \frac{\pi}{2}\cos \pi x\right]\cos \pi x$$

□

Integrating (3.16) gives us

(3.30) $$-\frac{2}{\pi}\sum_{n=1}^{\infty} \frac{1}{2n+1} \log\left(1+\frac{1}{n}\right) = -\int_0^1 \psi(x) \sin \pi x \, dx + \frac{2}{\pi}(\gamma + \log 2\pi)$$

Kölbig [20] obtained an equivalent result

$$\int_0^1 \psi(x) \sin \pi x \, dx = -\frac{2}{\pi}\left[\gamma + \log(2\pi) + 2\sum_{n=1}^{\infty} \frac{\log n}{4n^2 - 1}\right]$$

Kölbig [20] states that the above infinite series "does not seem to be expressible in terms of well-known functions". The equivalence is demonstrated below:

We have
$$2\sum_{n=1}^{\infty} \frac{\log n}{4n^2 - 1} = \sum_{n=1}^{\infty} \left[\frac{1}{2n-1} - \frac{1}{2n+1}\right] \log n$$

and consider the finite sum

$$\sum_{n=1}^{N} \frac{\log n}{2n-1} = \sum_{m=0}^{N-1} \frac{\log(m+1)}{2m+1}$$

$$= \sum_{n=1}^{N} \frac{\log(n+1)}{2n+1} - \frac{\log(N+1)}{2N+1}$$

Hence we have

$$\sum_{n=1}^{N} \left[ \frac{1}{2n-1} - \frac{1}{2n+1} \right] \log n = \sum_{n=1}^{N} \left[ \frac{\log(n+1)}{2n+1} - \frac{\log n}{2n+1} \right] - \frac{\log(N+1)}{2N+1}$$

Therefore, as $N \to \infty$ we see that

$$2\sum_{n=1}^{\infty} \frac{\log n}{4n^2-1} = \sum_{n=1}^{\infty} \frac{1}{2n+1} \log \frac{n+1}{n}$$

$$= \sum_{n=1}^{\infty} \frac{1}{2n+1} \log\left(1 + \frac{1}{n}\right)$$

$\square$

In 2006 Sondow and Hadjicostas [26] defined the generalised Euler constant function $\gamma(x)$ as

(3.31) $$\gamma(x) = \sum_{n=1}^{\infty} x^{n-1} \left[ \frac{1}{n} - \log\left(1 + \frac{1}{n}\right) \right]$$

where we see that $\gamma(1) = \gamma$ and, using (3.17), we have $\gamma(-1) = \log(4/\pi)$.

It is easily seen that

$$x^2 \gamma(x^2) = \sum_{n=1}^{\infty} x^{2n} \left[ \frac{1}{n} - \log\left(1 + \frac{1}{n}\right) \right]$$

and with $x \to e^{i\pi x}$ we have

$$e^{2i\pi x} \gamma(e^{2i\pi x}) = \sum_{n=1}^{\infty} \frac{e^{2ni\pi x}}{n} - \sum_{n=1}^{\infty} \log\left(1 + \frac{1}{n}\right) \cos 2n\pi x - i \sum_{n=1}^{\infty} \log\left(1 + \frac{1}{n}\right) \sin 2n\pi x$$

It should be noted that, by a very different method, Sondow and Hadjicostas [26] derived several formulae for a specific set of values of $\gamma(e^{i\pi x})$, namely $\gamma(e^{i\pi p/q})$ where $p$ and $q$ are relatively prime positive integers. One of their formulae is shown below.

$$\gamma(\omega) + \frac{\log(1-\omega)}{\omega} = \sum_{n=1}^{2q} \omega^{n-1} \log \frac{\Gamma\left(\frac{k+1}{2q}\right)}{\Gamma\left(\frac{k}{2q}\right)}$$

where $\omega = e^{i\pi p/q} \neq 1$.

We also have [26]

$$\gamma(x) = \int_0^1 \frac{1-y+\log y}{(1-xy)\log y} dy$$

□

As noted below, there is a connection between the series $\sum_{n=1}^{\infty} \frac{1}{n} \frac{1}{e^{2\pi nx}-1}$ and the Stieltjes constants.

Coffey [12] recently showed that

$$\gamma_1\left(\frac{3}{4}\right) - \gamma_1\left(\frac{1}{4}\right) = \pi\left[\frac{\pi}{3} + \gamma + 4\sum_{n=1}^{\infty} \frac{1}{n} \frac{1}{e^{2\pi n}-1}\right]$$

and using the known result

$$\gamma_1\left(\frac{3}{4}\right) - \gamma_1\left(\frac{1}{4}\right) = \pi\left[\gamma + 4\log 2 + 3\log \pi - 4\log \Gamma\left(\frac{1}{4}\right)\right]$$

we see that

$$4\sum_{n=1}^{\infty} \frac{1}{n} \frac{1}{e^{2\pi n}-1} = 4\log 2 + 3\log \pi - 4\log \Gamma\left(\frac{1}{4}\right) - \frac{\pi}{3}$$

The Euler reflection formula

$$\Gamma\left(\frac{1}{4}\right)\Gamma\left(\frac{3}{4}\right) = \pi\sqrt{2}$$

then gives us

$$\sum_{n=1}^{\infty} \frac{1}{n} \frac{1}{e^{2\pi n}-1} = \frac{1}{4}\log(4/\pi) + \log \Gamma\left(\frac{1}{4}\right) - \frac{\pi}{12}$$

which was previously recorded in Ramanujan's Notebooks as reported by Berndt [4, p.281].

□

Landau [22, p.181] showed that for $0 < x < \frac{1}{2}$

$$f\left(x + \frac{1}{2}\right) = f(2x) - f(x) - \log 2 \log(2\sin 2\pi x)$$

where $f(x)$ is defined as $f(x) = \sum_{n=1}^{\infty} \frac{\log n}{n} \cos 2n\pi x$. His proof proceeded as follows.

$$f\left(x+\frac{1}{2}\right) = \sum_{n=1}^{\infty} \frac{\log n}{n} \cos(2n\pi x + n\pi)$$

$$= \sum_{n=1}^{\infty} (-1)^n \frac{\log n}{n} \cos 2n\pi x$$

Subject to the necessary convergence conditions being satisfied we have

$$\sum_{n=1}^{\infty} (-1)^n a_n = -\sum_{n=1}^{\infty} a_n + 2\sum_{n=1}^{\infty} a_{2n}$$

and thus we have

$$\sum_{n=1}^{\infty} (-1)^n \frac{\log n}{n} \cos 2n\pi x = -\sum_{n=1}^{\infty} \frac{\log n}{n} \cos 2n\pi x + 2\sum_{n=1}^{\infty} \frac{\log(2n)}{2n} \cos 4n\pi x$$

$$= -f(x) + \log 2 \sum_{n=1}^{\infty} \frac{\cos 4n\pi x}{n} + \sum_{n=1}^{\infty} \frac{\log n}{n} \cos 4n\pi x$$

For $0 < x < \tfrac{1}{2}$ this becomes

$$= -f(x) - \log 2 \log(2\sin 2\pi x) + f(2x)$$

Therefore, using (3.22)

$$f(x) = \sum_{n=1}^{\infty} \frac{\log n}{n} \cos 2n\pi x = \frac{1}{2}\left[\varsigma''(0,x) + \varsigma''(0,1-x)\right] + \left[\gamma + \log(2\pi)\right]\log(2\sin \pi x)$$

we have the functional equation for $0 < x < \tfrac{1}{2}$

$$\varsigma''\left(0, x+\frac{1}{2}\right) + \varsigma''\left(0, \frac{1}{2}-x\right) = \varsigma''(0, 2x) + \varsigma''(0, 1-2x)$$

$$-\left[\varsigma''(0,x) + \varsigma''(0,1-x)\right] - 2\log 2 \log(2\sin 2\pi x)$$

Using (2.2) we obtain the following functional equation involving the first Stieltjes constants for $0 < x < 1/2$

$$\gamma_1\left(x+\frac{1}{2}\right) - \gamma_1\left(1-\left(x+\frac{1}{2}\right)\right) = 2\left[\gamma_1(2x) - \gamma_1(1-2x)\right] - \left[\gamma_1(x) - \gamma_1(1-x)\right] - 2\pi \log 2 \cdot \cot 2\pi x$$

## 4. A generalisation of the Briggs representation of the Stieltjes constants

Using the Poisson summation formula, we showed in [16] that

$$(4.1) \qquad \varsigma(s,x) = \frac{1}{2}\frac{1}{x^s} + \frac{x^{1-s}}{s-1} + 2\sum_{n=1}^{\infty}\int_0^{\infty}\frac{\cos 2\pi nt}{(x+t)^s}\,dt$$

and using $\frac{\partial}{\partial x}\varsigma(s,x) = -s\varsigma(s+1,x)$ we see that

$$\frac{\partial}{\partial x}\varsigma(s,x) = -\frac{1}{2}\frac{s}{x^{s+1}} - x^{-s} - 2s\sum_{n=1}^{\infty}\int_0^{\infty}\frac{\cos 2\pi nt}{(x+t)^{s+1}}\,dt$$

Differentiation results in

$$\frac{\partial^{m+1}}{\partial s^{m+1}}\frac{\partial}{\partial x}\varsigma(s,x) = -\frac{1}{2}s\frac{(-1)^{m+1}}{x^{s+1}}\log^{m+1}x - \frac{(-1)^m}{2}(m+1)\frac{1}{x^{s+1}}\log^m x - \frac{1}{x^s}(-1)^m\log^m x$$

$$-2(-1)^{m+1}s\sum_{n=1}^{\infty}\int_0^{\infty}\frac{\log^{m+1}(x+t)\cos 2\pi nt}{(x+t)^{s+1}}\,dt - 2(m+1)(-1)^m\sum_{n=1}^{\infty}\int_0^{\infty}\frac{\log^m(x+t)\cos 2\pi nt}{(x+t)^{s+1}}\,dt$$

and using (2.2)

$$\left.\frac{\partial^{m+1}}{\partial s^{m+1}}\frac{\partial}{\partial x}\varsigma(s,x)\right|_{s=0} = (-1)^{m+1}(m+1)\gamma_m(x)$$

we obtain for $m \geq 0$

$$(4.2) \qquad \gamma_m(x) = \frac{1}{2}\frac{\log^m x}{x} - \frac{\log^{m+1}x}{m+1} + 2\sum_{n=1}^{\infty}\int_0^{\infty}\frac{\log^m(x+t)}{x+t}\cos 2\pi nt\,dt$$

and this is a somewhat shorter exposition than the one shown in [16]. This formula was previously determined in an equivalent form by Zhang and Williams [30, Eq.(6.1)] in 1994.

We have for $m=0$

$$(4.3) \qquad \gamma_0(x) = -\psi(x) = \frac{1}{2x} - \log x + 2\sum_{n=1}^{\infty}\int_0^{\infty}\frac{\cos 2\pi nt\,dt}{x+t}$$

Whittaker & Watson [29, p.261] posed the following question: Prove that for all values of $x$ except negative real values we have

$$(4.4) \qquad \log\Gamma(x) = \frac{1}{2}\log(2\pi) + \left(x - \frac{1}{2}\right)\log x - x + \frac{1}{\pi}\sum_{n=1}^{\infty}\int_0^{\infty}\frac{\sin(2n\pi t)}{n(x+t)}\,dt$$

and this result was attributed by Stieltjes to Bourguet. Equation (5.4) may also be derived using the Euler-Maclaurin summation formula (see for example Knopp's book [21, p.530]).

We note that (4.3) may be obtained by differentiating (4.4).

We have for $m \geq 1$

$$(4.5) \qquad \gamma_m = 2 \sum_{n=1}^{\infty} \int_0^{\infty} \frac{\log^m(1+x)}{1+x} \cos 2\pi nx \, dx$$

and with $x \to 1+x$ we obtain

$$(4.6) \qquad \gamma_m = 2 \sum_{n=1}^{\infty} \int_1^{\infty} \frac{\log^m x}{x} \cos 2\pi nx \, dx$$

as originally shown by Briggs [7].

## 5. Another series for the generalised Stieltjes constants

From the treatise by Srivastava and Choi [27, p.146] we find a different expression for the Hurwitz zeta function

$$(5.1) \qquad \varsigma(s,x) = \frac{x^{1-s}}{s-1} - \sum_{n=1}^{\infty} (-1)^n \frac{\Gamma(s+n)}{(n+1)\Gamma(n+1)\Gamma(s)} \varsigma(s+n,x)$$

Using $\frac{\partial}{\partial x} \varsigma(s,x) = -s\varsigma(s+1,x)$ we see that

$$\frac{\partial}{\partial x} \varsigma(s,x) = -x^{-s} + s \sum_{n=1}^{\infty} (-1)^n \frac{\Gamma(s+n+1)}{(n+1)\Gamma(n+1)\Gamma(s+1)} \varsigma(s+n+1,x)$$

$$:= -x^{-s} + sG(s,x)$$

where

$$G(s,x) := \sum_{n=1}^{\infty} (-1)^n \frac{g_n(s,x)}{(n+1)\Gamma(n+1)}$$

and

$$g_n(s,x) := \frac{\Gamma(s+n+1)}{\Gamma(s+1)} \varsigma(s+n+1,x)$$

Differentiation gives us

$$\frac{\partial^{m+1}}{\partial s^{m+1}} \frac{\partial}{\partial x} \varsigma(s,x) = -(-1)^{m+1} x^{-s} \log^{m+1} x + s \frac{\partial^{m+1}}{\partial s^{m+1}} G(s,x) + (m+1) \frac{\partial^m}{\partial s^m} G(s,x)$$

and with $s=0$ we have

$$\frac{\partial^{m+1}}{\partial s^{m+1}}\frac{\partial}{\partial x}\varsigma(s,x)\bigg|_{s=0} = -(-1)^{m+1}\log^{m+1}x + (m+1)\frac{\partial^m}{\partial s^m}G(s,x)\bigg|_{s=0}$$

We define $h(s)$ as

$$h(s) = \frac{\Gamma(s+n+1)}{\Gamma(s+1)}$$

and note that

$$h'(s) = \frac{\Gamma(s+n+1)}{\Gamma(s+1)}[\psi(s+n+1) - \psi(s+1)]$$

Therefore, the higher derivatives of $h(s)$ may be represented in terms of the (exponential) complete Bell polynomials

$$h^{(k)}(s) = \frac{\Gamma(s+n+1)}{\Gamma(s+1)} Y_k\left(\psi(s+n+1)-\psi(s+1), \psi'(s+n+1)-\psi'(s+1),...,\psi^{(k-1)}(s+n+1)-\psi^{(k-1)}(s+1)\right)$$

The (exponential) complete Bell polynomials may be defined by $Y_0 = 1$ and for $r \geq 1$

$$Y_r(x_1,...,x_r) = \sum_{\pi(r)} \frac{r!}{k_1! \, k_2! ... k_r!} \left(\frac{x_1}{1!}\right)^{k_1} \left(\frac{x_2}{2!}\right)^{k_2} ... \left(\frac{x_r}{r!}\right)^{k_r}$$

where the sum is taken over all partitions $\pi(r)$ of $r$, i.e. over all sets of integers $k_j$ such that

$$k_1 + 2k_2 + 3k_3 + \cdots + rk_r = r$$

Further details of the (exponential) complete Bell polynomials are contained in [9] and [13].

With $s = 0$ we have

$$h^{(k)}(0) = \Gamma(n+1) Y_k\left(\psi(n+1)-\psi(1), \psi'(n+1)-\psi'(1),...,\psi^{(k-1)}(n+1)-\psi^{(k-1)}(1)\right)$$

$$\frac{\partial^m}{\partial s^m} g_n(s,x)\bigg|_{s=0} = \sum_{k=0}^{m} \binom{m}{k} h^{(k)}(0) \varsigma^{(m-k)}(n+1,x)$$

$$\frac{\partial^{m+1}}{\partial s^{m+1}}\frac{\partial}{\partial x}\varsigma(s,x)\bigg|_{s=0} = -(-1)^{m+1}\log^{m+1}x + (m+1)\sum_{n=1}^{\infty}(-1)^n \frac{1}{n+1}\sum_{k=0}^{m}\binom{m}{k} h^{(k)}(0)\varsigma^{(m-k)}(n+1,x)$$

Using (2.2) we obtain

$$\gamma_m(x) = -\frac{1}{m+1}\log^{m+1}x + (-1)^{m+1}\sum_{n=1}^{\infty}(-1)^n \frac{1}{n+1}\sum_{k=0}^{m}\binom{m}{k} h^{(k)}(0)\varsigma^{(m-k)}(n+1,x)$$

We have the well-known formula [27, p.14]

$$\psi(n+t) - \psi(t) = H_n^{(1)}(t)$$

where $H_n^{(m)}(t)$ is the generalised harmonic number function defined by

$$H_n^{(m)}(x) = \sum_{k=0}^{n-1} \frac{1}{(k+t)^m}$$

and we note that

$$H_n^{(m)}(1) = H_n^{(m)} = \sum_{k=1}^{n} \frac{1}{k^m}$$

We have the derivatives

$$\frac{d^r}{dt^r}[\psi(n+t) - \psi(t)] = (-1)^r r! H_n^{(r+1)}(t)$$

$$\left.\frac{d^r}{dt^r}[\psi(n+t) - \psi(t)]\right|_{t=1} = (-1)^r r! H_n^{(r+1)}$$

Hence, we obtain

$$h^{(k)}(0) = \Gamma(n+1) Y_k\left(H_n^{(1)}, -1! H_n^{(2)}, \ldots, (-1)^k k! H_n^{(k)}\right)$$

and, for convenience, we denote

$$\mathbf{Y}_k(H_n^{(k)}) := Y_k\left(H_n^{(1)}, -1! H_n^{(2)}, \ldots, (-1)^k k! H_n^{(k)}\right)$$

We therefore have

(5.2) $$\gamma_m(x) = -\frac{1}{m+1}\log^{m+1} x + (-1)^{m+1}\sum_{n=1}^{\infty}(-1)^n \frac{1}{n+1}\sum_{k=0}^{m}\binom{m}{k}\mathbf{Y}_k(H_n^{(k)})\varsigma^{(m-k)}(n+1,x)$$

and

(5.3) $$\gamma_m = (-1)^{m+1}\sum_{n=1}^{\infty}(-1)^n \frac{1}{n+1}\sum_{k=0}^{m}\binom{m}{k}\mathbf{Y}_k(H_n^{(k)})\varsigma^{(m-k)}(n+1)$$

This generalises the following results previously obtained by Coffey [12]

$$\gamma_0(x) = -\log x - \sum_{n=1}^{\infty}\frac{(-1)^n}{n+1}\varsigma(n+1,x)$$

$$\gamma_1(x) = -\frac{1}{2}\log^2 x + \sum_{n=1}^{\infty}(-1)^n \frac{\varsigma'(n+1,x)}{n+1} + \sum_{n=1}^{\infty}(-1)^n \frac{H_n \varsigma(n+1,x)}{n+1}$$

## 9. Open access to our own work

This paper contains references to various other papers and, rather surprisingly, most of them are currently freely available on the internet. Surely now is the time that <u>all</u> of <u>our</u> work should be freely accessible by <u>all</u>. The mathematics community should lead the way on this by publishing <u>everything</u> on arXiv, or in an equivalent open access repository. We think it, we write it, so why hide it? You know it makes sense.

Wessex House,
Devizes Road,
Upavon,
Pewsey,
Wiltshire SN9 6DL